\documentclass[12pt]{amsproc}
\usepackage{hyperref}
\usepackage{geometry}
\usepackage{graphicx}
\usepackage{amsthm}
\usepackage{amssymb}

\theoremstyle{plain}
\newtheorem{thm}{Theorem}[section]
\newtheorem{cor}[thm]{Corollary}
\newtheorem{lem}[thm]{Lemma}
\newtheorem{prop}[thm]{Proposition}
\newtheorem{defn}[thm]{Definition}
\newtheorem{exa}[thm]{Example}
\newtheorem{rem}[thm]{Remark}

\begin{document}

\title{GENERALIZATIONS OF GRADED $S$-PRIMARY IDEALS}

\author{Tamem \textsc{Al-Shorman}}
\address{Department of Mathematics and Statistics, Jordan University of Science and Technology, Irbid, Jordan}
\email{alshorman91@yahoo.com}

\author{Malik \textsc{Bataineh}}
\address{Department of Mathematics and Statistics, Jordan University of Science and Technology, Irbid, Jordan}
\email{msbataineh@just.edu.jo}

\author{Rashid \textsc{Abu-Dawwas}}
\address{Department of Mathematics, Yarmouk University, Irbid, Jordan}
\email{rrashid@yu.edu.jo}



\begin{abstract}
The goal of this article is to present the graded weakly $S$-primary ideals and $g$-weakly $S$-primary ideals which are extensions of graded weakly primary ideals. Let $R$ be a commutative graded ring, $S\subseteq h(R)$ and $P$ be a graded ideal of $R$. We state $P$ is a graded weakly $S$-primary ideal of $R$ if there exists $s\in S$ such that for all $x,y \in h(R)$, if $0\neq xy\in P$, then $sx\in P$ or $sy\in Grad(P)$ (the graded radical of $P$). Several properties and characteristics of graded weakly $S$-primary ideals as well as graded $g$-weakly $S$-primary ideals are investigated.
\end{abstract}

\maketitle

\section{Introduction}

 Throughout this article, $G$ will be an abelian group  with identity $e$ and $R$ be a commutative ring with nonzero unity $1$ element. $R$ is called a $G$-graded ring if $ R= \bigoplus\limits_{g \in G} R_g$   with the property $R_gR_h\subseteq R_{gh}$ for all $g,h \in G$, where $R_g$ is an additive subgroup of $R$ for all $g\in G$. The elements of $R_g$ are called homogeneous of degree $g$. If $x\in R$, then $x$ can be written uniquely as $\sum\limits_{g\in G} x_g$, where $x_g$ is the component of $x$ in $R_g$. The set of all homogeneous elements of $R$ is $h(R)= \bigcup\limits_{g\in G} R_g$. Let $P$ be an ideal of a $G$-graded ring $R$. Then $P$ is called a graded ideal if $P=\bigoplus\limits_{g\in G}P_g$, i.e, for $x\in P$ and  $x=\sum\limits_{g\in G} x_g$ where $x_g \in P_g$ for all $g\in G$. An ideal of a $G$-graded ring is not necessary graded ideal (see \cite{abu2019graded}). The concept of graded prime ideals and its generalizations have an indispensable role in commutative $G$-graded rings.

Let $P$ be a proper graded ideal of $R$. Then the graded radical of $P$ is denoted by $Grad(P)$ and it is defined as written below:
\begin{center}
  {\small $Grad(P)=\Big{\{} x= \sum\limits_{g\in G} x_g \in R$ : for all $g\in G$, there exists $n_g\in \mathbb{N}$ such that ${x_g}^{n_g} \in P \Big{\}}$}.
\end{center}
Note that $Grad(P)$ is always a graded ideal of $R$ (see \cite{refai2000graded}).

Let $P$ be a graded ideal of $R$ and $P\neq R$, then $P$ is said to be graded primary ideal of $R$, if $x,y\in h(R)$ with $xy\in P$ then $x\in P$ or $y \in Grad(P)$ (see \cite{refai2004graded}). The concept of a graded primary ideal can be generalized in a many ways, for example, Bataineh in \cite{bataineh2011generalizations} and Atani in \cite{atani2005graded}. A proper graded ideal $P$ of $R$ is claimed to be a graded weakly primary ideal if whenever $0\neq xy \in P$ for some $x,y \in h(R)$, then either $x\in P$ or $y \in Grad(P)$.

Let $S\subseteq R$ be a multiplicative set and $P$ be an ideal of $R$ such that $P\cap S= \emptyset$. For an ideal $P$ of $R$, the radical of $P$ denoted by $rad(P)$ is the ideal $\{ x\in R : x^n\in P$ for some positive integer $n$ $\}$ of $R$. In \cite{massaoud2021s}, $P$ is said to be an $S$-primary ideal of $R$ if there exists $s\in S$ such that for all $x,y \in R$, if $xy\in P$, then $sx\in P$ or $sy\in rad(P)$. $P$ is said to be a weakly $S$-primary ideal of $R$ if there exists $s\in S$ such that for all $x,y \in R$, if $0\neq xy \in P$, then $sx\in P$ or $sy \in rad(P)$ (see \cite{celikeli2021weakly}).

Let $R$ be a $G$-graded ring, $S\subseteq h(R)$ be a multiplicative set and $P$ be a graded ideal of $R$ such that $P\cap S = \emptyset$. Recently, Azzh in \cite{alshehry2021graded}, defined the concept of graded $S$-primary ideals which is another generalization of graded primary ideals. A proper graded ideal $P$ of $R$ is said to be a graded $S$-primary ideal of $R$ if there exists $s\in S$ such that for all $x,y \in h(R)$, if $xy \in P$, then $sx\in P$ or $sy \in Grad(P)$. Also, we study the relationship between graded weakly $S$-primary ideals and other graded ideals. We will classify all $G$-graded rings for which they proper graded primary ideal is a graded weakly $S$-primary ideal. In Section Two, motivated by \cite{celikeli2021weakly}, we introduce the concept of graded weakly $S$-primary ideal. We say that $P$ is a graded weakly $S$-primary ideal of $R$ if there exists $s\in S$ such that for all $x,y \in h(R)$, if $0\neq xy \in P$, then $sx\in P$ or $sy \in Grad(P)$. In Proposition \ref{prop1}, if $S$ consists of units of $h(R)$, then graded weakly primary ideals and graded weakly $S$-primary ideals coincide. We show that every graded ideal is a graded weakly $S$-primary ideal in a $G$-graded valuation domain. We investigate some basic properties of graded weakly $S$-primary ideals.
\\

\begin{center}
Weakly $S$-primary $\rightarrow$ Graded weakly $S$-primary.

\vspace{0.25cm}

The inverse is not always true (Example \ref{ex1}).

\vspace{0.25cm}
Graded $S$-primary $\rightarrow$ Graded weakly $S$-primary.

\vspace{0.25cm}

The inverse is not always true (Example \ref{ex2}).

\vspace{0.25cm}
Graded weakly primary $\rightarrow$ Graded weakly $S$-primary.

\vspace{0.25cm}

The inverse is not always true (Example \ref{ex3}).
\end{center}

In Section Three, we introduce the concept of $g$-$S$-primary ideal. We say that $P$ is a $g$-$S$-primary ideal of $R$ if there exists $s\in S\subseteq R_e$ such that for all $x,y \in R_g$, if $ xy \in P$, then $sx\in P$ or $sy \in Grad(P)$, and we introduce the concept of $g$-weakly $S$-primary ideal and give some of its basic properties. We say that $P$ is a $g$-weakly $S$-primary ideal of $R$ if there exists $s\in S\subseteq R_e$ such that for all $x,y \in R_g$, if $0\neq xy \in P$, then $sx\in P$ or $sy \in Grad(P)$. Finally, we show that every graded ideal is a $g$-weakly $S$-primary ideal in $G$-graded valuation domain.

\section{Graded Weakly $S$-Primary Ideals}

In this section, we introduce the concept of graded weakly $S$-primary ideals and show that graded weakly $S$-primary ideals have a lot of similar features to these of graded weakly primary ideals.

\begin{defn}\label{def1}
Let $R$ be a $G$-graded ring, $S\subseteq h(R)$ be a multiplicative set and $P$ be a graded ideal of $R$ such that $P\cap S =\emptyset$. We say that $P$ is a graded weakly $S$-primary ideal of $R$, if there exists $s\in S$ such that for all $x,y\ \in\ h(R)$, if $0\ \neq \ xy \ \in\ P$,\ then $sx\ \in\ P$ or $sy\ \in\ Grad(P)$.
\end{defn}

 Obviously, every weakly $S$-primary ideal is graded weakly $S$-primary ideal. However, the converse is not always true.

\begin{exa}\label{ex1}\cite{alshehry2021graded}
Let $ R = \mathbb{Z}$[i] and $G = \mathbb{Z}_2$. So, R is G-graded ring. Let $I=5R$ be graded ideal of R. To see more details to show I is graded prime ideal of R  (see \cite{alshehry2021graded}, Example 1).
\\
Now, let $P=10R$ be a graded ideal of R and $S =\{ 2^n : n \in \mathbb{Z}^+ \} \subseteq h(R)$. It was shown that P is graded S-primary ideal of R, then P is graded weakly S-primary ideal of R. However, P is not weakly S- primary ideal of R. Now, because (7-i),(7+i) $\in$ R such that $0 \neq (7-i)(7+i) \in P$, but $(s(7-i))^n \notin P $ and $(s(7+i))^n \notin P $ for some s $\in $ S and $n \in \mathbb{Z}^+$.
 \end{exa}
 
Obviously, every graded $S$-primary ideal is graded weakly $S$-primary ideal. However, the opposite is not always true.

 \begin{exa}\label{ex2}
Let $ R=\mathbb{Z}_{12}[i]$ and $G=\mathbb{Z}_2$. Then R is G-graded ring. Let $P=\{ 0 \}$ be graded ideal of R and  $S=\{ 1,3,9\} \subseteq h(R)$. P is graded weakly S-primary ideal of R where $P\cap S =\emptyset $, but P is not graded S-primary ideal of R because 3,4 $\in$ h(R) and 3.4 $\in$ P, but $(3s)^n \notin P $ and $(4s)^n \notin P$ for some $s\in S$.
 \end{exa}

 Moreover, it is evident that every graded weakly primary ideal and it is disjoint to S is a graded weakly S-primary ideal. However, the opposite is not always true.

 \begin{exa}\label{ex3}\cite{alshehry2021graded}
 Consider $R=\mathbb{Z}[X]$ and $G=\mathbb{Z}$. Then R is G-graded ring by $R_j =\mathbb{Z} X^j$ for $j \geq 0$ and $R_j = \{ 0 \}$ otherwise. Let $P=9XR$ be a graded ideal of R and  $S= \{ 3^n :  n \in \mathbb{Z}^+  \} \subseteq h(R)$. Note that $P\cap S =\emptyset $. Therefor, P is graded S-primary ideal of R, and hence P is a graded S-primary ideal of R. Since P is a graded S-primary ideal of R, so P is graded weakly S-primary ideal of R. However, P is not graded weakly primary ideal of R, because 18,X$\in$ h(R) and $0\neq 18X\in P$, but $(18)^n \notin P$ and $X^n \notin P$, $n \in \mathbb{Z}^+$.
 \end{exa}

 \begin{prop}\label{prop1}
 Let R be a G-graded ring and S \  $\subseteq$ \ h(R) be a multiplicative set. If S consists of units of h(R), then graded weakly primary ideals and graded weakly S-primary ideals are coincide.
 \\
\textbf{Proof.} By contradiction, let P be a graded weakly primary ideal of R and $P'$ be a graded weakly S-primary ideal of R where $P \neq P'$.
 \\
 Since P is graded weakly primary ideal of R, so for all x,y $\in$ h(R) with $0\neq xy \in P$, then $x\in P$ or $y \in Grad(P)$.
 \\
 Since $P'$ is graded weakly S-primary ideal of R, so for all x,y $\in$ h(R) and there exists  s $\in$ S with $0\neq xy \in P'$, so sx $\in$ $P'$ or sy $\in$ Grad($P'$).
 \\
 Case1: if sx $\in$ $P'$, since s $\in$ S, so there exists  $s^{-1} \in S$ such that $s^{-1} s = 1_R \in h(R)$. So $s^{-1} s x = 1_R x =x \in P'$ since $1_R \in S $.
 \\
 case2: if sy $\in$ Grad($P'$), since s $\in$ S, so there exists  $s^{-1} \in S$ such that $s^{-1} s = 1_R \in h(R)$. So $s^{-1 } s y = 1_R y =y \in Grad(P')$ since $1_R \in S $.
 \\
 Thus, by case1 and case2 we have $P'=P$ for all x,y $\in$ h(R). Hence a graded weakly primary ideal coincide of graded weakly S-primary ideal.
 \end{prop}

Note that in R a G-graded ring and I is a graded ideal of R, then $I \subseteq Grad(I)$ and  Grad(Grad(I))=Grad(I) (see \cite{refai2000graded}).

 \begin{prop}\label{prop2}
 Let R be a G-graded ring, S$\subseteq$ h(R) be a multiplicative set and P is a graded ideal of R such that $P \cap S= \emptyset$. Then P is a graded weakly S-primary ideal of R if and only if Grad(P) is a graded weakly S-prime ideal of R.
\\
 \textbf{Proof.} ($\Rightarrow$) Since $P \cap S = \emptyset$, so Grad(P) $\cap$ $S = \emptyset$. Let for all x,y $\in$ h(R) with $0\neq xy \in Grad(P)$. So $0 \neq (xy)^n = x^ny^n \in P$ for some $n\in \mathbb{Z}^+$, and so there exists  $s \in S$ such that $sx^n \in P $ or $sy^n \in Grad(P)$, this means that $sx \in Grad(P)$ or $sy\in Grad(Grad(P))= Grad(P)$. Therefore, Grad(P) is a graded weakly S-prime ideal of R.
 \\
 ($\Leftarrow$) Assume that Grad(P) is graded weakly S-prime ideal of R and for all x,y $\in$ h(R) with $0\neq xy  \in P $, then $0\neq xy \in Grad(P)$ and there exists $s \in S$ such that $sx\in Grad(P)$ or $sy \in Grad(P)$. So, $sx^n \in P$ for some $n\in \mathbb{Z}^+$, and so $sx \in P$. Thus P is graded weakly S-primary ideal of R.
\end{prop}

  Note that $S^{-1}R$ is a G-graded ring with $S^{-1}R=\{ \frac{r}{s}: r\in R , s\in S\}$. If P is a graded ideal of R and $P \cap S = \emptyset$,  then $S^{-1}P \neq S^{-1}R$ and $S^{-1}Grad(P)=Grad(S^{-1}P)$ and $S^{-1}P$ is a graded ideal of $S^{-1}R$ (see \cite{nastasescu2004methods} and  \cite{atiyah1969introduction}).

\begin{thm}\label{thm1}
   Let R be a G-graded ring and S$\subseteq$ h(R) be a multiplicative set. If P is a graded weakly S-primary ideal of R, then $S^{-1}P$ is a graded weakly primary ideal of $S^{-1}R$.
   \\
   \textbf{Proof.} Since $S \cap P = \emptyset$, we have $S^{-1}P \neq S^{-1}R$. Let $0 \neq \frac{x}{s_1} \frac{y}{s_2} \in S^{-1}P$ for all $x,y \in h(R)$ and $s_1,s_2 \in S$. Then $ \frac{x}{s_1} \frac{y}{s_2} = \frac{z}{s_3}$ for some $z \in P$ and $s_3 \in S$. So, there  is $t \in S$ such that $0 \neq ts_2xy=ts_1s_2z \in P$. As P is a graded weakly S-primary ideal of R, then there exists  $s\in S$ such that $sts_3 \in P$ or $sxy \in Grad(P)$. Thus $0\neq sxy \in Grad(P)$ and $sts_3 \notin P$ because $S \cap P = \emptyset$. Hence, $0 \neq s^2x \in Grad(P)$ or $sy \in Grad(P)$ because Grad(P) is a graded S-prime ideal of R, and so $sx \in Grad(P)$ or $sy \in Grad(P)$. So,  $\frac{x}{s_1}= \frac{sx}{ss_1} \in S^{-1} Grad(P)$ or $\frac{y}{s_2}= \frac{sy}{ss_2} \in S^{-1} Grad(P)$. But $S^{-1}Grad(P) = Grad(S^{-1}P)$, so $\frac{x}{s_1} \in Grad(S^{-1}P)$, then ${(\frac{x}{s_1})}^n \in S^{-1}P$ for some positive integer n, take $n=1$, so $\frac{x}{s_1} \in S^{-1}P$ or $\frac{y}{s_2} \in S^{-1}Grad(P)= Grad(S^{-1}P)$. Therefore, $S^{-1}P$ is a graded weakly primary ideal of $S^{-1}R$.
\end{thm}

Next proposition generalize the notion of graded S-primary ideal. First, we needed the following lemma.

 \begin{lem}\label{lem1}
  Let R be a G-graded ring and S$\subseteq$ h(R) be a multiplicative set. If P is a graded weakly S-primary ideal of R such that $P \cap S= \emptyset$. Then $Grad((P:s)) = Grad((P:s^n))$ for some $s \in S$ and positive integer n.
  \\
  \textbf{Proof.} To see the details of the proof, see \cite{alshehry2021graded} in Proposition 2.
 \end{lem}

 \begin{prop}\label{prop3}
    Let R be a G-graded ring, S$\subseteq$ h(R) be a multiplicative set and P is a graded ideal of R such that $P \cap S= \emptyset$. Then P is a graded weakly S-primary ideal of R if and only if (P:s) is a graded weakly primary ideal of R for some $s \in S$.
  \\
    \textbf{Proof.} Similar to Proposition 2 in \cite{alshehry2021graded}
 \end{prop}

On the contrary of result in Proposition\ref{prop3}, the requirement that S be made up of regular elements was required. The following example demonstrates

\begin{exa}\label{ex4}
As in Example \ref{ex3}, P=$\{0\}$ is graded weakly S-primary ideal of R. As may be seen, $s=1 \in S$ with (P:s)=P is graded primary ideal of R. However, $3\in S$ is a zero divisor since 3.4=0 and $4 \in R$.
\end{exa}

\begin{prop}\label{prop4}
Let R be a G-graded ring, S$\subseteq$ h(R) be a multiplicative set consisting of regular elements and P is graded ideal of R such that $P \cap S= \emptyset$. Then P is graded weakly S-primary ideal of R if and only if $S^{-1}P$ is graded weakly primary ideal of $S^{-1}R$ and there is $s\in S$ such that $(P:t) \subseteq (P:s)$ for all $t\in S$.
\\
\textbf{Proof.} $(\Rightarrow)$ Assume that P is graded weakly S-primary ideal of R and for all $x, y \in h(R)$  there exists $s\in S$ with $0 \neq xy\in P$. Then there is $sx \in P$ or $sy \in Grad(P)$. Thus, by Theorem \ref{thm1}, $S^{-1}P$ is a graded weakly primary ideal of $S^{-1}R$. Let $t\in S$ and $0\neq x \in (P:t)$. So $0\neq tx \in P$, and so $sx \in P $ or $st \in Grad(P)$. Since $P \cap S= \emptyset$ so, $st \notin Grad(P)$ which is implies that $sx\in P$, hence $x\in (P:s)$.
\\
($\Leftarrow$) Let $0\neq xy \in P$ for all $x,y \in h(R)$. Then $0 \neq \frac{x}{1}\frac{y}{1}\in S^{-1}P$ for some, and then $\frac{x}{1}\in S^{-1}P$ or $\frac{y}{1} \in Grad(S^{-1}P)$. If $\frac{x}{1} \in S^{-1}P$, then $\frac{x}{1}=\frac{w_1}{t}$ for some $w_1\in P $, and then $xt = w_1 \in P$, which implies that $x\in (P:t) \subseteq (P:s)$. Hence, $sx\in P$. If  $\frac{y}{1} \in Grad(S^{-1}P)$, so, $\frac{y^n}{1^n}=\frac{y^n}{1} \in S^{-1}P$ then $\frac{y^n}{1}=\frac{w_2}{t}$ for some $w_2 \in P$,and then $y^nt = w_2 \in P$, which implies that $ y^n \in (P:t) \subseteq (P:s)$. Then $y^n \in Grad((P:s))=Grad((p:s^n))$ by Lemma \ref{lem1}, and then $s^ny^n =(sy)^n \in p$. Thus $sy \in Grad(P)$. Therefore, P is a graded weakly S-primary ideal of R.
\end{prop}

Note that in R be a G-graded ring and $S\subseteq h(R)$ be a multiplicative set of R, if P and I are a graded ideal of R, then $S^{-1}(P:I)=(S^{-1}P:S^{-1}I)$ (see \cite{atiyah1969introduction}).

\begin{thm}\label{thm2}
   Let R be a G-graded ring, S$\subseteq$ h(R) be a multiplicative set and P is a graded ideal of R such that $P \cap S= \emptyset$. Then P is a graded weakly S-primary ideal of R if and only if $S^{-1}P$ is a graded weakly primary ideal of $S^{-1}R$ and $S^{-1}P=(P:s)$ for some $s\in S$.
\\
\textbf{Proof.} $(\Rightarrow)$ Assume that P is a graded weakly S-primary ideal of R. So there exists $s\in S$ and for all $x, y \in h(R)$ with $0 \neq xy\in P$, we have $sx \in P$ or $sy \in Grad(P)$. Then by Theorem \ref{thm1}, $S^{-1}P$ is graded weakly primary ideal of $S^{-1}R$. Now, let $a\in (P:s)$, then $sa\in P$, and then $a= \frac{sa}{a} \in S^{-1}P$, and hence $(P:s) \subseteq S^{-1}P $. Let $b \in S^{-1}P $ and $b=\frac{w}{t}$ for some $w\in P$ and $s\in S$, then $0 \neq bt=w \in P$. So, $sb \in P$ or $st \in Grad(P)$. Since $P \cap S = \emptyset$, then $st \notin Grad(P)$, and then $sb \in P$, we get  $b \in (P:s)$. Therefore, $S^{-1}P=(P:s)$ for some $s\in S$.
\\
($\Leftarrow$) Since $S^{-1}P$ is a graded weakly primary ideal of R and $S^{-1}P =(P:s)$ so, (P:s) is a graded weakly primary ideal of R. Then by Proposition \ref{prop3} P is a graded weakly S-primary ideal of R.
\end{thm}

\begin{rem}\label{rem1}
In Example \ref{ex3}, we show that P is a graded weakly S-primary ideal of R. So, $S^{-1}P$ is a graded weakly primary ideal of $S^{-1}R$, whilst P is not graded primary ideal of R.
\end{rem}

\begin{thm}\label{thm3}
    Let R be a G-graded ring, S$\subseteq$ h(R) be a multiplicative set and P is a graded ideal of R such that $P \cap S= \emptyset$. Then the following statements are equivalent:
   \\
   $(i)$ P is a graded weakly S-primary ideal of R.
   \\
   $(ii)$ (P:s) is a graded weakly primary ideal of R for some $s \in S$.
   \\
   $(iii)$ $S^{-1}P$ is a graded weakly primary ideal of $S^{-1}R$ and $S^{-1}P=(P:s)$ for some $s\in S$.
   \\
   \textbf{Proof.} Follows by Proposition \ref{prop3} and Theorem \ref{thm2}.
\end{thm}

Note that in R a G-graded ring and I,J are graded ideals of R, then $I\cap J$ is a graded ideal of R (see \cite{farzalipour2012union}).

\begin{prop}\label{prop5}
Let R be a G-graded ring, S$\subseteq$ h(R) be a multiplicative set and P is a graded weakly S-primary ideal of R such that $P \cap S= \emptyset$. If I is graded ideal of R such that $I\cap S \neq \emptyset$, then $P\cap I$ is graded weakly S-primary ideals of R.
\\
\textbf{Proof.} Since P and I are graded ideals of R, then $P\cap I$ is a graded ideal of R. Clearly, $(P \cap I) \cap S = \emptyset$. Let $t \in I \cap S$ and $0 \neq xy \in P \cap I$ for all $x,y \in h(R)$. Then $0\neq xy \in P$, and then $sx\in P$ or $sy \in Grad(P)$ for some $s\in S$. This means $w=st \in S$. If $sx \in P$, so,  $wx \in P\cap I$, or if $sy \in Grad(P)$, so, $sy^n \in P$. Then $wy^n \in P\cap I$, and then $ y^n \in (P\cap I : w)=(P\cap I : w^n)$. Hence $w^ny^n= (wy)^n \in P\cap I$, then $wy \in Grad(P\cap I)$ for some $w \in S$. Therefore, $P\cap I$ is a graded weakly S-primary ideal of R.
\end{prop}

\begin{rem}\label{rem2}
Let $S_1 \subseteq S_2$ be multiplicative subset of h(R) and P be a graded ideal of R such that  $P \cap S_2 = \emptyset$. If P is a graded weakly $S_1$-primary ideal of R, then P is graded weakly $S_2$-primary ideal of R. However, the converse is not always true.
\\
\textbf{Proof.} We show in Example \ref{ex3} that $P=9XR$ is graded weakly $S_2$-primary ideal of R where $S_2 = \{ 3^n$: n is non-negative integer$\}$, however, $S_1=\{1\} \subseteq S_2$ with P is not a graded weakly $S_1$-primary ideal of R, since 18,X$\in$ h(R) with $0\neq 18X\in P$, but $(18)^n \notin P$ and $X^n \notin P$, n positive integer.
\end{rem}

\begin{prop}\label{prop6}
Let R be a G-graded ring and $S_1,S_2\subseteq$ h(R) be a multiplicative set with $S_1 \subseteq S_2$. For any $s\in S_2$, there is $t \in S_2$ satisfying $st\in S_1$. If P is a graded weakly $S_2$-primary ideal of R, then P is a graded weakly $S_1$-primary ideal of R.
\\
\textbf{Proof.} Since P is graded weakly $S_2$-primary ideal of R, then there exists $s\in S_2$ such that if  $0\neq xy \in P$ for all $x,y \in h(R)$, so $sx \in P$ or $sy \in Grad(P)$. By assumption $w=ts\in S_1$ for some $t\in S_2$, then $wx\in P$ or $wy \in Grad(P)$. Consequently, P is a graded weakly $S_1$-primary ideal of R.
\end{prop}

Let $S\subseteq h(R)$ be a multiplicative set and $S^* =\{ r\in h(R) : \frac{r}{1}$ is unit in $S^{-1}R\}$ denotes the saturation of S. Note that, $S^* \subseteq h(R)$ is a multiplicative set containing S (see \cite{saber2021graded}).

\begin{prop}\label{prop7}
Let R be a G-graded ring, S$\subseteq$ h(R) be a multiplicative set and P is a graded ideal of R such that $P \cap S= \emptyset$. Then P is a graded weakly S-primary ideal of R if and only if P is a graded weakly $S^*$-primary ideal of R.
\\
\textbf{Proof.} ($\Rightarrow$) Clearly, since $S\subseteq S^*$ and $P \cap S^* = \emptyset$. By Remark \ref{rem2}, P is graded weakly $S^*$-primary ideal of R.
\\
($\Leftarrow$) Just we proved that for any $s\in S^*$, there is $t \in S^*$ such that $st\in S$. Let $s\in S^*$, then $\frac{s}{1}\frac{a}{s_1}=1$ for some $s_1 \in S$ and $a \in h(R)$. This means $s_2as=s_2s_1\in S$ for some $s_2 \in S$. Now, take $t=s_2a$. Then, we have $t \in S^*$ with $st \in S$. Therefore, by putting $S_1=S$ and $S_2=S^*$. We deduce the result from Proposition \ref{prop6}.
\end{prop}

In the next proposition we show that the intersection of graded weakly S-primary ideals is also graded weakly S-primary ideal. Firstly, we needed the following lemma.

\begin{lem}\label{lem2}
  Let R be a G-graded ring and S$\subseteq$ h(R) be a multiplicative set. If $P_i$  for all  $i=\{1,...,n\}$ are a graded weakly S-primary ideals of R such that $P_i \cap S= \emptyset$ for all $i=\{1,...,n\}$. Then $\bigcap\limits_{i=1}^{n} Grad(P_i) = Grad(\bigcap\limits_{i=1}^{n} P_i )  $.
  \\
  \textbf{Proof.} Let $x\in \bigcap\limits_{i=1}^{n} Grad(P_i)$, then $x \in Grad(P_i)$  for all $i=\{1,...,n\}$, and then $x^k \in P_i$ for some positive integer k. So, $x^k \in \bigcap\limits_{i=1}^{n} P_i$. Hence, $x \in Grad(\bigcap\limits_{i=1}^{n}P_i)$. conversely, let $y \in Grad(\bigcap\limits_{i=1}^{n} P_i)$, then $y^k \in \bigcap\limits_{i=1}^{n} P_i$ for some positive integer k, and then $y^k \in P_i$ for all $i=\{1,...,n\}$. So, $y \in Grad(P_i)$ for all $i=\{1,..,n\}$. Thus $y \in \bigcap\limits_{i=1}^{n} Grad(P_i)$.
\end{lem}

\begin{prop}\label{prop8}
Let R be a G-graded ring, S$\subseteq$ h(R) be a multiplicative set and $P_i$  for all  $i=\{1,...,n\}$ are  graded ideals of R such that $P_i \cap S= \emptyset$. If $P_i$ is a graded weakly S-primary ideal of R for each $i$ with $Grad(P_i)=Grad(P_j)$ for all $i,j \in \{1,...,n\}$, then $\bigcap\limits_{i=1}^{n} P_i $ is a graded weakly S-primary ideal of R.
\\
\textbf{Proof.} Since $P_i$ is a graded weakly S-primary ideal of R, there exists $s_i \in S$ and for all $x,y \in h(R)$ with $0\neq xy \in P_i$, we have either $s_ix \in P_i$ or $s_iy \in Grad(P_i)$. let $s= \prod\limits_{i=1}^{n} s_i$. Then $s\in S$. Assume that $x,y \in h(R)$ and $0 \neq xy \in \bigcap\limits_{i=1}^{n} P_i $, then $0\neq xy \in P_i $ for all $i=\{1,...,n\}$. Since $P_i$ is a graded weakly S-primary ideals of R, then $sx\in P_i$ or $sy \in Grad(P_i)$ for each $i$. Thus $sx \in \bigcap\limits_{i=1}^{n} P_i $ or $sy \in \bigcap\limits_{i=1}^{n} Grad(P_i) = Grad(\bigcap\limits_{i=1}^{n} P_i)$. Therefore, $\bigcap\limits_{i=1}^{n} P_i$ is a graded weakly S-primary ideal of R.
\end{prop}

Note that let R be a G-graded ring. Then the graded radical of R is denoted by Grad(R) and equal R (see \cite{refai2000graded}).

\begin{lem}\label{lem3}
  Let $R=R_1 \times R_2$ where each $R_i$ is a G-graded ring for all $i=\{1,2\}$. Then the following hold:
  \\
  $(i)$ $P_1$ is a graded ideal of $R_1$, if and only if $Grad(P_1 \times R_2) = Grad(P_1) \times R_2$.
  \\
  $(ii)$ $P_2$ is a graded ideal of $R_2$, if and only if $Grad(R_1 \times P_2) = R1 \times Grad(P_2)$.
  \\
  $(iii)$ $P_1$ is a graded ideal of $R_1$ and $P_2$ is a graded ideal of $R_2$, if and only if $Grad(P_1 \times P_2) = Grad(P_1) \times Grad(P_2)$.
  \\
  \textbf{Proof.} $(i)$ Let $(x,y) \in Grad(P_1 \times R_2)$, then $(x,y)^n \in P_1 \times R_2 $, and then $(x,y)^n= (x^n,y^n) \in P_1 \times R_2$. So, $x^n \in P_1$ and $y^n \in R_2$, thus $x \in Grad(P_1)$ and $y \in Grad(R_2)=R_2$. Therefore, $(x,y) \in Grad(P_1) \times R_2$. Conversely, let $(s,t) \in Grad(P_1)\times R_2$, then $s \in Grad(P_1)$ and $t \in R_2 = Grad(R_2)$, and then $s^n \in P_1$ and $t^n\in R_2$. So, $(s^n,t^n)=(s,t)^n \in P_1 \times R_2$, hence $(s,t) \in Grad(P_1 \times R_2)$. Therefore, $Grad(P_1 \times R_2) = Grad(P_1) \times R_2$.
  \\
  $(ii)$ Similar to  case $(i)$.
  \\
  $(iii)$  Let $(x,y) \in Grad(P_1 \times P_2)$, then $(x,y)^n \in P_1 \times P_2 $, and then $(x,y)^n= (s^n,t^n) \in P_1 \times P_2$. So, $s^n \in P_1$ and $t^n \in P_2$, thus $s \in Grad(P_1)$ and $t \in Grad(P_2)$. Therefore, $(x,y) \in Grad(P_1) \times Grad(P_2)$. Conversely, let $(s,t) \in Grad(P_1)\times Grad(P_2)$, then $s \in Grad(P_1)$ and $t \in Grad(P_2)$, and then $s^n \in P_1$ and $t^n\in P_2$. So, $(s^n,t^n)=(s,t)^n \in P_1 \times P_2$, hence $(s,t) \in Grad(P_1 \times P_2)$. Therefore, $Grad(P_1 \times P_2) = Grad(P_1) \times Grad(P_2)$.
\end{lem}

\begin{thm}\label{thm4}
   Let $R=R_1 \times R_2$ where each $R_i$ is a G-graded ring for all $i=\{ 1,2 \}$ and $S =S_1 \times S_2 \subseteq h(R) $ be a multiplicative set such that $P_1 \cap S_1 = \emptyset $ and $P_2 \cap S_2 = \emptyset $. Then the following hold:
   \\
   $(i)$ $P_1$ is a graded weakly $S_1$-primary ideal of $R_1$ if and only if  $P_1 \times R_2$ is a graded weakly S-primary ideal of R.
   \\
   $(ii)$ $P_2$ is a graded weakly $S_2$-primary ideal of $R_2$ if and only if  $R_1 \times P_2$ is graded weakly S-primary ideal of R.
   \\
   $(iii)$ $P_1$ is graded weakly $S_1$-primary ideal of $R_1$ and $P_2$ is graded weakly $S_2$-primary ideal of $R_2$ if and only if  $P_1 \times P_2$ is graded weakly S-primary ideal of R.
   \\
   \textbf{Proof.}$(i)$ Let $0 \neq (a,b)(c,d) =(ac,bd) \in P_1 \times R_2$ where $(a,b),(c,d)\in  h(R)$, since $P_1$ is a graded weakly $S_1$-primary ideal of $R_1$, then either $sa \in P_1$ or $sc \in Grad(P_1)$ for some $s\in S_1$. Let $t \in S_2$, so $t \in R_2$ and $tb \in R_2  $ and $td \in R_2$. It follows that either $(sa,tb)=(s,t)(a,b) \in P_1 \times R_2$ or $(sc,td)=(s,t)(c,d) \in Grad(P_1) \times R_2 = Grad(P_1 \times R_2)$ for some $(s,t) \in S$. Therefore,  $P_1 \times R_2$ is a graded weakly S-primary ideal of R. Conversely, take $(a,b), (c,d) \in h(R)$ and $0\neq (a,b)(c,d) \in P_1 \times R_2$ then $0\neq ac \in P_1$. So, $(s,t)(a,b) \in P_1\times R_2$ or $(s,t)(c,d)\in Grad(P_1\times R_2)= Grad(P_1) \times Grad(R_2)$ for some $(s,t)\in S = S_1 \times S_2$. Then $sa \in P_1$ or $sc \in Grad(P_1)$. Hence $P_1$ is a graded weakly $S_1$-primary ideal of $R_1$.
   \\
   $(ii)$ It is similar to $(i)$.
   \\
   $(iii)$ Let $0 \neq (a,b)(c,d) =(ac,bd) \in P_1 \times P_2$ where $(a,b),(c,d)\in  h(R)$, since $P_1$ is a graded weakly $S_1$-primary ideal of $R_1$, then either $sa \in P_1$ or $sc \in Grad(P_1)$ for some $s\in S_1$ and since $P_2$ is graded weakly $S_2$-primary ideal of $R_2$, so either $tb \in P_2$ or $td \in Grad(P2)$ for some $t\in S_2$. It follows that either $(sa,tb)=(s,t)(a,b) \in P_1 \times P_2$ or $(sc,td)=(s,t)(c,d) \in Grad(P_1) \times Grad(P_2) = Grad(P_1 \times P_2)$ for some $(s,t) \in S$. Therefore, $P_1 \times P_2$ is a graded weakly S-primary ideal of R. Conversely, take $(a,b), (c,d) \in h(R)$ and $0\neq (a,b)(c,d) \in P_1 \times P_2$ then $0\neq ac \in P_1$ and $0\neq bd \in P_2$. So, $(s,t)(a,b) \in P_1\times P_2$ or $(s,t)(c,d)\in Grad(P_1\times P_2)= Grad(P_1) \times Grad(P_2)$ for some $(s,t)\in S = S_1 \times S_2$. Then $sa \in P_1$ or $sc \in Grad(P_1)$ thus $P_1$ is a graded weakly $S_1$-primary ideal of $R_1$, and then $tb \in P_2$ or $td \in Grad(P_1)$ hence $P_2$ is a graded weakly $S_2$-primary ideal of $R_2$.
\end{thm}

\begin{thm}\label{thm5}
   Let $R=R_1 \times R_2$ where each $R_i$ is a G-graded ring for all $i=\{ 1,2 \}$ and $S =S_1 \times S_2 \subseteq h(R) $ be a multiplicative set such that $P_1 \cap S_1 = \emptyset $ and $P_2 \cap S_2 = \emptyset $. Then the following statements are equivalent:
   \\
   $(i)$ $P_1 \times P_2$ is graded weakly S-primary ideal of R.
   \\
   $(ii)$ $P_1=R_1$ and $P_2$ is a graded weakly $S_2$-primary ideal of $R_2$ or $P_2=R_2$ and $P_1$ is a graded weakly $S_1$-primary ideal of $R_1$ or $P_1$ is a graded weakly $S_1$-primary ideal of $R_1$ and $P_2$ is a graded weakly $S_2$-primary ideal of $R_2$.
   \\
   \textbf{Proof.} It follows from Theorem \ref{thm4}.
\end{thm}

\begin{prop}\label{prop9}
Let R be a G-graded ring, S$\subseteq$ h(R) be a multiplicative set and P is graded ideal of R such that $P \cap S= \emptyset$. Then P is graded weakly S-primary ideal of R if and only if there exists $s\in S$ and for all graded ideals I,J of R, if $IJ \subseteq P$, then $sI \subseteq P$ or $ sJ \subseteq P$.
\\
\textbf{Proof.} ($\Rightarrow$) Assume that P is graded weakly S-primary ideal of R. So there exists $s\in S$ and for all $x,y \in h(R)$ with $0\neq xy \in P$, there is $sx\in P$ or $sy \in Grad(P)$. Let I,J be graded ideals of R with $IJ\subseteq P$. Suppose that $sI \not\subseteq P$. Then $sx\notin P$ for any $x \in I$. Let $y \in J$, then $xy \in IJ \subseteq P$, and then $sy \in Grad(P)$ for any $y\in J$. So, $(sy)^n \in P$ for some positive integer n. Take n=1 so $sy\in P$, then $sJ\subseteq P$.
\\
($\Leftarrow$) Let for all $x,y\in h(R)$ with $0\neq xy \in P$. Then I=Rx and J=Ry are graded ideals of R with $IJ\subseteq P$, so  by assumption. $sI \subseteq P$ or $ sJ \subseteq P$, and hence $sx \in P$ or $sy \in P$ so $sy \in Grad(P)$. Thus, P is graded weakly S-primary ideal of R.
\end{prop}

Next corollary can be deduced from Proposition \ref{prop9} using induction:

\begin{cor}\label{coro1}
Let R be a G-graded ring, S$\subseteq$ h(R) be a multiplicative set and P is a graded ideal of R such that $P \cap S= \emptyset$. Then P is a graded weakly S-primary ideal of R if and only if there exists $s\in S$ and for all graded ideals $I_1,I_2,...,I_n$ of R, if $I_1I_2...I_n \subseteq P$, then $sI_i \subseteq P$ for some $1 \leq i \leq n$.
\end{cor}

Recall that if $R_1$ and $R_2$ are G-graded rings. A ring homomorphism $f:R_1 \rightarrow R_2$ is said to be graded ring homomorphism if $f(R_1)\subseteq R_2$ (see \cite{nastasescu2004methods}). Let $f:R_1 \rightarrow R_2$ be graded ring homomorphism and I, J are two graded ideal of $R_1$, $R_2$ respectively, then $f(I)$ , $f^{-1}(J)$ be graded ring of $R_2$, $R_1$ respectively (see \cite{refai2020generalizations}).

\begin{prop}\label{prop10}
Let $f:R_1 \rightarrow R_2$ be a graded ring homomorphism, S$\subseteq$ h($R_1$) be a multiplicative set and P,
I are two graded ideal of $R_1$, $R_2$ respectively such that $P\cap S=\emptyset$ and $I\cap f(S)=\emptyset$. Then the following hold:
\\
$(i)$ If  $f$  is a graded epimorphism and P is a graded weakly S-primary ideal of $R_1$ containing  Ker(f), then $f(P)$ is a graded weakly $f(S)$-primary ideal of $R_2$.
\\
$(ii)$ If $f$ is a graded monomorphism and I is a graded weakly $f(S)$-primary ideal of $R_2$, then $f^{-1}(I)$ is a graded weakly S-primary ideal of $R_1$.
\\
\textbf{Proof.} $(i)$ Let $f(P)$ be a graded ideal of $R_2$ and $r\in f(P)\cap f(S)$. Then $r=f(p)=f(s)$ when $p\in P$ and $s\in S$, and then $s-p \in Ker(f)\subseteq P$ this means $s\in P$, this is contradiction. Hence $f(P) \cap f(S)=\emptyset$. Now, let $0\neq xy \in f(P)$ for all $x,y \in h(R_2)$. Then there some $a,b \in h(R_1)$ such that f(a)=x and f(b)=y, so $0\neq f(ab)=xy \in f(P)$. Since $Ker(f)\subseteq P$, we obtain $0\neq ab \in P$, and then $sa\in P$ or $sb\in Grad(P)$ for some $s\in S$. This means that $f(s)x\in f(P)$ or $f(s)y \in Grad(f(P))$. Hence, $f(P)$ is graded weakly $f(S)$-primary ideal of $R_2$.
\\
$(ii)$ Let $f^{-1}(I)$ be graded ideal of $R_1$. So, $f^{-1}(I)\cap S=\emptyset$. Now, let $x,y \in h(R_1) $ and $0\neq xy \in f^{-1}(I)$. Since $Ker(f)=\{0\}$, then $0\neq f(xy)=f(x)f(y) \in I$, and then $f(s)f(x)=f(sx)\in I$ or $f(s)f(y)=f(sy)\in Grad(I)$ for any $s\in S$. Hence, $sx\in f^{-1}(I)$ or $sy \in Grad(f^{-1}(I))$. As a result, we can deduce that $f^{-1}(I)$ is graded weakly S-primary ideal of $R_1$.
\end{prop}

In R is G-graded ring and S$\subseteq$ h(R). If P and I  are  graded ideals of R, then P/I is graded ideal of $R/I$ and $\Bar{S}=\{s+I:s\in S\}$ is a multiplicative subset of $R/I$ (see \cite{saber2020graded} and \cite{celikeli2021weakly}).

\begin{cor}\label{coro2}
Let R be a G-graded ring, S$\subseteq$ h(R) be a multiplicative set and P, I are two graded ideal of R with $I\subseteq P$. Then the following hold:
\\
$(i)$ If P is a graded weakly S-primary ideal of R, then $P/I$ is a graded weakly $\Bar{S}$-primary ideal of $R/I$.
\\
$(ii)$ If P is a graded weakly S-primary ideal of R, then $P\cap R_e$ is a graded weakly S-primary ideal of $R_e$.
\\
$(iii)$ If $P/I$ is a graded weakly $\Bar{S}$-primary ideal of $R/I$ and I is a graded S-primary ideal of R, then P is a graded S-primary ideal of R.
\\
$(iv)$ If $P/I$ is a graded weakly $\Bar{S}$-primary ideal of $R/I$ and I is a graded S-primary ideal of R, then P is a graded weakly S-primary weakly ideal of R.
\\
\textbf{Proof.} $(i)$ Let $f:R\rightarrow R/I$ be graded epimorphism by $f(r)=r+I$. So by the result follows by Proposition \ref{prop10} $(i)$.
\\
$(ii)$ Let $f:R_e\rightarrow R$ be graded monomorphism by $f(r)=r$. So by the result follows by Proposition \ref{prop10} $(ii)$.
\\
$(iii)$ Let $P/I$ is a graded weakly $\Bar{S}$-primary ideal of $R/I$ and I is a graded S-primary ideal of R. Let $ xy \in P$ for all $x,y \in h(R)$. If $xy \in I$, then $sx \in I \subseteq P$ or $sy \in Grad(I) \subseteq Grad(P)$ for any $s\in S$. If $xy \in I$, then $I\neq (x+I)(y+I) \in P/I$ this implies that $(t+I)(x+I)\in P/I$ or $(t+I)(y+I)\in Grad(P/I)=Grad(P)/I$ for some $(t+I)\in \Bar{S}$. Thus $tx \in P$ or $ts \in Grad(P)$. Consequently, P is a graded S-primary ideal of R.
\\
$(iv)$ Is similar to $(iii)$.
\end{cor}

In R a G-graded ring and I,J and P are graded ideals of R, then $I + J$ is a graded ideal of R and we say P is graded maximal ideal of R if $P\neq R$ and there is no graded ideal I of R such that $P \subset I \subset R$  (see \cite{farzalipour2012union} and \cite{refai2000graded}).

\begin{prop}\label{prop11}
Let R be a G-graded ring and S$\subseteq$ h(R) be a multiplicative set. If P,I are two graded weakly S-primary ideals of R with $I\subseteq P$ such that $(P+I)\cap S =\emptyset$, then P+I is a graded weakly S-primary ideal of R.
\\
\textbf{Proof.} Let $\Bar{S}=\{s+I:s\in S\}$ and $\Bar{S_1}=\{s+(P\cap I):s\in S\}$ be  multiplicative subset of $R/(P\cap I)$ and $R/I$, respectively. So, $P/(P\cap I)$ is graded weakly $\Bar{S_1}$-primary ideal of $R/(P\cap I)$ by Corollary \ref{coro2}. Since $\Bar{S_1}\subseteq \Bar{S}$, then $P/(P\cap I)$ is also graded weakly $\Bar{S}$-primary ideal of $R/(P\cap I)$ by Remark \ref{rem2}. Thus, $(P+I)/I$  is a graded weakly $\Bar{S}$-primary ideal of $R/ I$ by isomorphic between $(P+I)/I$ and $P/(P\cap I)$. As a result of Corollary \ref{coro2} that is graded weakly S-primary ideal of R.
\end{prop}

\begin{lem}\label{lem4}
Let R be a G-graded ring, S$\subseteq$ h(R) be a multiplicative set and I be a graded ideal of R such that $I\cap S=\emptyset$. Then there is a graded ideal P of R which is maximal with $I\subseteq P$ and $P\cap S =\emptyset$. Moreover, P is a graded primary ideal.
\\
\textbf{Proof.}  To see the details of the proof, see \cite{saber2021graded} in Lemma 3.10. But has been switched from graded prime ideal to graded primary ideal.
\end{lem}

If every nonzero homogeneous element of a G-graded ring R is unit, it is said to be a G-graded field. Similarly, if R does not have a homogeneous zero divisor, it is said to be a G-graded domain (see \cite{saber2020graded}).

\begin{prop}\label{prop12}
 Let R be a G-graded ring and S$\subseteq$ h(R) be a multiplicative set. Then the following statements are equivalent:
\\
$(i)$ $\{0\}$ is the only graded weakly S-primary ideal of R.
\\
$(ii)$ $\{0\}$ is the only graded S-primary ideal of R.
\\
$(iii)$ R is a G-graded domain and $S^{-1}R$ is a G-graded field.
\\
\textbf{Proof.} $(i) \Rightarrow (ii)$: Let $P=\{0\}$ be graded S-primary ideal of R, so P is graded weakly S-primary ideal of R. Hence, $\{0\}$ is only graded S-primary ideal of R.
\\
$(ii) \Rightarrow (iii)$: By Lemma \ref{lem4}, there is a primary ideal P of R that has been graded and $P\cap S= \emptyset$. Thus P is graded S-primary ideal of R. So $P=\{0\}$, then R is G-graded domain. Take $0\neq x\in h(R)$ and $s\in S$. We proved that $\frac{x}{s}$ is unit in $S^{-1}R$ . If $x\in S$, so we have the outcome that we want. Suppose that  $x\notin S$. If $R_x\cap S = \emptyset$, so by Lemma \ref{lem4}, there is a primary ideal P of R that has been graded such that $R_x\subseteq P =\{0\}$, a contradiction. Then $R_x\cap S \neq \emptyset$. Let $t\in R_x \cap S$, so $t\in S$ and t=rx for any $r \in R$. We have, $\frac{sr}{t} \in S^{-1}R$ and $\frac{x}{s}\frac{sr}{t}=\frac{xsr}{st}=\frac{st}{st}= 1$. Then $\frac{x}{s}$ is unit in $S^{-1}R$, and hence $S^{-1}R$ is a G-graded field.
\\
$(iii) \Rightarrow (i)$: Let P is a nonzero graded weakly S-primary ideal of R. Let $0\neq p \in P$. Since $S^{-1}R$ is G-graded field, then there exists $0\neq x \in R$ and $s\in S$ such that $\frac{px}{s}=1$. Since R is a G-graded domain, we deduce that $px=s \in P \cap S$, a contradiction. Hence, $\{0\}$ is only graded weakly S-primary ideal of R.
\end{prop}

\begin{thm}\label{thm6}
   Let R be a G-graded ring and S$\subseteq$ h(R) be a multiplicative set. Then every graded weakly S-primary ideal of R is graded primary if and only if R is a G-graded domain and every graded S-primary ideal of R is graded primary.
   \\
   \textbf{Proof.} ($\Rightarrow$) Assume that every graded weakly S-primary ideal of R is graded primary. Because $\{0\}$ is graded weakly S-primary ideal of R, then $\{0\}$ is graded primary ideal of R, and then R is G-graded domain. In addition, every graded S-primary ideal of R is graded weakly S-primary ideal of R,and thus is graded primary by hypothesis.
\\
($\Leftarrow$) Since R is G-graded domain, every graded weakly S-primary ideal of R is graded S-primary ideal of R, and hence is graded primary by hypothesis.
\end{thm}

Let R be a G-graded ring and $S\subseteq h(R)$ be a multiplicative set. Then R is called a G-graded valuation domain if for each nonzero homogeneous elements $a,b \in h(R)$ either $a|b$ or $b|a$, and then the set of graded ideals of R is totally ordered under inclusion (see \cite{uregen2019graded}). Let I be a graded ideal of R, then Grad(I) is the intersection of all graded prime ideal of R containing I (see \cite{refai2000graded}). The intersection of graded prime ideal of R is also graded prime ideal of R and every graded prime ideal of R that is disjoint with S is a graded weakly S-prime ideal of R (see \cite{saber2021graded}).

\begin{thm}\label{thm7}
   Let R be a G-graded ring and S$\subseteq$ h(R) be a multiplicative set. If R is a G-graded valuation domain , then every graded ideal P of R such that $P\cap S=\emptyset$ is a graded weakly S-primary ideal  of R.
   \\
   \textbf{Proof.} Let P be a graded ideal of R, then $Grad(P)= \bigcap\limits_{P\subseteq I} I$, where I denotes any graded prime ideal containing P. Since R is a G-graded valuation domain all graded prime ideals which contains P is linearly ordered, So Grad(P)=$I_\alpha$ for some graded prime ideal $I_\alpha$, and so Grad(P) is a graded prime ideal of R thus Grad(P) is a graded weakly S-prime ideal of R. By Proposition \ref{prop2} P is a graded weakly S-primary ideal of R.
\end{thm}

\section{$g$-WEAKLY $S$-PRIMARY IDEALS}
In this section, we introduce the concept of $g$-weakly S-primary ideals and show that $g$-weakly $S$-primary ideals have a lot of similar features to these of graded weakly $S$-primary ideals.

\begin{defn}\label{def2}
A graded ideal $P$ of a $G$-graded ring $R$ is $g$-$S$-primary ideal if there exists $s\in S \subseteq R_e$ when $P\cap S =\emptyset$ such that for all $a,b \in R_g$ when $g\in G$, if $ab \in P$, then $sa\in P$ or $sb\in Grad(P)$.
\end{defn}

\begin{defn}\label{def3}
A graded ideal $P$ of a $G$-graded ring $R$ is $g$-weakly $S$-primary ideal if there exists $s\in S \subseteq R_e$ when $P\cap S =\emptyset$ such that for all $a,b \in R_g$ when $g\in G$, if $0\neq ab \in P$, then $sa\in P$ or $sb\in Grad(P)$.
\end{defn}

\begin{lem}\label{lem5}
If P is a graded weakly S-primary ideal in G-graded ideal R, then P is a graded g-weakly S-primary ideal of R.
\\
\textbf{Proof.} For $g \in G$ and $s\in S \subseteq R_e$, assume that $0\neq xy \in P$ when $x,y \in R_g \subseteq h(R)$, so either $sx\in P$ or $sy \in Grad(P)$ since P is a graded weakly S-primary ideal of R. It follows that P is a graded g-weakly S-primary ideal of R.
\end{lem}

Clearly, any graded g-S-primary ideal is also a graded g-weakly S-primary ideal, but this is not always the case, as shown in the example below.

\begin{exa}\label{ex5}
Let $ R=\mathbf{Z}_{12}[i]$ and $G=\mathbf{Z}_2$, then R is G-graded ring by $R_0=\mathbf{Z}_{12}$. Consider a graded ideal $P=\{ 0 \}$ of R and let multiplicative subset $S=\{ 1,3,9\}$ of $R_0$. Note that $P\cap S =\emptyset $. P is a graded weakly S-primary ideal of R, so P is graded g-weakly S-primary ideal for all $g\in G$ by Lemma \ref{lem5}. On the other hand, P is not a graded 0-S-primary ideal of R since $3,4 \in R_0$ with $3.4 \in P$, $(3s)^n \notin P $ and $(4s)^n \notin P$ for each $s\in S$.
\end{exa}

\begin{prop}\label{prop13}
 Let R be a G-graded ring, then P is a graded g-weakly S-primary ideal of R if and only if Grad(P) is a graded g-weakly S-prime ideal of R.
 \\
 \textbf{Proof.}($\Rightarrow$) $Grad(P) \cap S =\emptyset$ because $p\cap S =\emptyset$. Let $x, y \in R_g$ when $g\in G$ such that $0\neq xy \in Grad(P)$. So  $0 \neq (xy)^n = x^ny^n \in P$ for some $n \in \mathbf{Z}^+$, and there exist $s\in S\subseteq R_e$ such that $sx^n \in P $ or $sy^n \in Grad(P)$, implying that $sx \in Grad(P)$ or $sy\in Grad(Grad(P))= Grad(P)$. Therefore, Grad(P) is a graded g-weakly S-prime ideal of R.
 \\
 ($\Leftarrow$) Suppose that Grad(P) is graded g-weakly S-prime ideal of R. Let $x,y \in R_g$ when $g\in G$ such that $0\neq xy  \in P $, so $0\neq xy \in Grad(P)$, and so there exists  $s \in S \subseteq R_e$ such that $sx\in Grad(P)$ or $sy \in Grad(P)$. So, $sx^n \in P$ for some $n\in \mathbf{Z}^+$, and so $sx \in P$. Thus P is graded g-weakly S-primary ideal of R.
\end{prop}

\begin{prop}\label{prop14}
Let P be a graded g-weakly S-primary ideal of G-graded ring R. If I is a graded ideal of R such that $I\cap S \neq \emptyset$ when $S\subseteq R_e$, then $P\cap I$ is a graded g-weakly S-primary ideals of R.
\\
\textbf{Proof.}Since P and I are graded ideals of R, then $P\cap I$ is a graded ideal of R. Clearly, $(P \cap I) \cap S = \emptyset$. Let $t \in I \cap S$ and $0 \neq xy \in P \cap I$ for some $x,y \in R_g $ when $g\in G$. Then $0\neq xy \in P$, and then $sx\in P$ or $sy \in Grad(P)$ for some $s\in S\subseteq R_e$. implying that $w=st \in S$. If $sx \in P$, so,  $wx \in P\cap I$, or if $sy \in Grad(P)$, so  $sy^n \in P$. Then $wy^n \in P\cap I$, and then $ y^n \in (P\cap I : w)=(P\cap I : w^n)$. Hence $w^ny^n= (wy)^n \in P\cap I$, then $wy \in Grad(P\cap I)$ for some $w \in S$. Therefore, $P\cap I$ is a graded g-weakly S-primary ideal of R.
\end{prop}

\begin{rem}\label{rem3}
Let P be a graded ideal of G-graded ring R and $S_1\subseteq S_2 \subseteq R_e$ be a multiplicative set of R such that  $P \cap S_2 = \emptyset$. If P is a graded g-weakly $S_1$-primary ideal of R, then P is graded g-weakly $S_2$-primary ideal of R. However, as shown in the following example, the opposite is not always true.
\end{rem}

\begin{exa}\label{ex6}
Let $R=\mathbf{Z}[X]$ be a G-graded ring when $G=\mathbf{Z}$ by $R_j =\mathbf{Z} X^j$ for $j \geq 0$ and $R_j = \{ 0 \}$ otherwise. Take the graded ideal $P=9XR$ of R and $S_2= \{ 3^n :  n \in \mathbf{Z}^+ \} \subseteq R_e$ is a multiplicative set of R. Then P is a graded g-weakly $S_2$-primary ideal of R. However, $S_1=\{1\} \subseteq S_2$ with P is not a graded g-weakly $S_1$-primary ideal of R, since $27,X \in R_g$ when $g\in G$ with $0\neq 27X\in P$, but $27^n \notin P$ and $X^n \notin P$ for some $n\in \mathbf{Z}^+$.
\end{exa}

\begin{prop}\label{prop15}
Let $S_1,S_2\subseteq R_e$ be a multiplicative sets of G-graded ideal R with $S_1 \subseteq S_2$. For any $s\in S_2$, there is an $t \in S_2$ satisfying $st\in S_1$. If P is a graded g-weakly $S_2$-primary ideal of R, then P is a graded g-weakly $S_1$-primary ideal of R.
\\
\textbf{Proof.} Since P is a graded g-weakly $S_2$-primary ideal of R, there exists $s\in S_2$ such that if  $0\neq xy \in P$ for some $x,y \in R_g$ when $g\in G$, then $sx \in P$ or $sy \in Grad(P)$. By assumption $w=ts\in S_1$ for some $t\in S_2$, then $wx\in P$ or $wy \in Grad(P)$. Consequently, P is a graded g-weakly $S_1$-primary ideal of R.
\end{prop}

\begin{prop} \label{prop16}
Let $S \in R_e$ be a multiplicative set of G-graded ring R and P be a graded ideal of R such that $P \cap S= \emptyset$. Then P is a graded g-weakly S-primary ideal of R if and only if P is a graded g-weakly $S^*$-primary ideal of R.
\\
\textbf{Proof.} ($\Rightarrow$) Is clearly, since $S\subseteq S^*$ and $P \cap S^* = \emptyset$. By Remark \ref{rem3}, P is a graded g-weakly $S^*$-primary ideal of R.
\\
($\Leftarrow$) Just we show that for any $s\in S^*$, there is $t \in S^*$ such that $st\in S$. Let $s\in S*$, then $\frac{s}{1}\frac{a}{s_1}=1$ for some $s_1 \in S$ and $a \in R_g$ when $g\in G$. This implies that $s_2as=s_2s_1\in S$ for some $s_2 \in S$. Now, take $t=s_2a$. Then, we have $t \in S^*$ with $st \in S$. Therefore, by putting $S_1=S$ and $S_2=S^*$. We deduce the result from Proposition \ref{prop15}.
\end{prop}

\begin{prop}\label{prop17}
Let $S \subseteq R_e$ be a multiplicative set of G-graded ring R and $P_i$  for all  $i=\{1,...,n\}$ be a graded ideals of R such that $P_i \cap S= \emptyset$. If $P_i$ is a graded g-weakly S-primary ideal of R for each $i$ with $Grad(P_i)=Grad(P_j)$ for all $i,j \in \{1,...,n\}$, then $\bigcap\limits_{i=1}^{n} P_i $ is a graded g-weakly S-primary ideal of R.
\\
\textbf{Proof.} Since $P_i$ is a graded g-weakly S-primary ideal of R, there exists $s_i \in S$ and for all $x,y \in R_g$ when $g\in G$ with $0\neq xy \in P_i$, we have either $s_ix \in P_i$ or $s_iy \in Grad(P_i)$. let $s= \prod\limits_{i=1}^{n} s_i$. Then $s\in S$. Assume that $x,y \in R_g$ when $g\in G$ and $0 \neq xy \in \bigcap\limits_{i=1}^{n} P_i $, then $0\neq xy \in P_i $ for all $i=\{1,...,n\}$. Since $P_i$ is a graded g-weakly S-primary ideals of R, then $sx\in P_i$ or $sy \in Grad(P_i)$ for each $i$. Thus $sx \in \bigcap\limits_{i=1}^{n} P_i $ or $sy \in \bigcap\limits_{i=1}^{n} Grad(P_i) = Grad(\bigcap\limits_{i=1}^{n} P_i)$. Therefore, $\bigcap\limits_{i=1}^{n} P_i$ is a graded g-weakly S-primary ideal of R.
\end{prop}

\begin{prop}\label{prop18}
If P be a graded g-weakly S-primary ideal of G-graded ring R which is not graded g-S-primary ideal of R, then $P_g^2 =\{0\}$.
\\
\textbf{Proof.} By contradiction, assume that $P_g^2 \neq \{0\}$ and P is a graded g-weakly S-primary ideal. Let $x,y \in R_g$ such that $xy \in P$. If $xy \neq 0$, then $sx\in P$ or $sy \in Grad(P)$ when $s\in S\subseteq R_e$. Suppose $xy=0$. If $xP_g\neq \{0\}$, there is $p\in P_g$ such that $xp\neq 0$, and so $0\neq xy = x(p+y)\in P$. Hence, $sx\in P$ or $s(p+y)\in Grad(P)$. As $p\in P \subseteq Grad(P)$, we have $sx\in P$ or $sy \in Grad(P)$. Similarly, if $yP_g\neq \{0\}$, we arrive to the same conclusion.. So, we can assume that $xP_g=\{0\}$ and $yP_g=\{0\}$. Since $P_g^2 \neq \{0\}$, there exists $p,q \in P_g \subseteq P \subseteq Grad(P)$ such that $pq\neq 0$. Thus, $0\neq pq =(x+p)(y+q) \in P$. Then $s(x+p)\in P$ or $s(y+q)\in Grad(P)$, and hence either $sx\in P$ or $sy \in P$. Consequently, we conclude that P is a graded g-S-primary ideal of R. So, $P_g^2=\{0\}$.
\end{prop}

\begin{cor}\label{coro3}
  Let P be a g-weakly primary ideal of R which is not g-primary ideal of R. Then $P_g^2=\{0\}$.
  \\
 \textbf{Proof.} By Proposition \ref{prop18} and take $S=\{1\}$.
\end{cor}

\begin{cor}\label{coro4}
 Let P be a graded g-weakly S-primary ideal of R which is not graded g-S-primary ideal of R. Then $P_g\subseteq Grad(\{0\})$. In particular, if R is reduced, then $P_g=\{0\}$.
\end{cor}

\begin{thm}\label{thm8}
Let $S\subseteq R_e$ be a multiplicative set of G-graded ring R and P be a graded ideal of R such that $P\cap S=\emptyset$. Then the following statements are equivalent:
\\
$(i)$ P is a graded g-weakly S-primary ideal of R.
\\
$(ii)$ For all $a\notin (P:_{R_g} s)$ when $s\in S$ then we have either $(P:_{R_g} a) \subseteq (P:_{R_g} s)$ or $(P:_{R_g} a)=(0:_{R_g} a)$.
\\
$(iii)$ For all graded ideal I and J of R and there exists $s\in S$, if $0\neq I_g J_g\subseteq P$, then $sI_g\subseteq P$ or $sJ_g\subseteq P$.
\\
\textbf{Proof.} Similar to Theorem 3.26 in \cite{saber2021graded}.
\end{thm}

\begin{prop}\label{prop19}
Let P be a graded e-weakly S-primary ideal of G-graded ring R that is not graded e-S-primary ideal of R. Then $sP_e(Grad(\{0\}))_e=\{0\}$ for some $s\in S$.
\\
\textbf{Proof.} By Theorem \ref{thm8} and similar to Propositin 3.27 in \cite{saber2021graded}.
\end{prop}

\begin{cor}\label{coro5}
  Let P be a graded e-weakly S-primary ideal of G-graded ring R that is not graded e-S-primary ideal of R. Then $P_g\subseteq Grad(\{0\})$ and $P_g(Grad(\{0\}))_e=\{0\}$.
  \\
  \textbf{Proof.} By Corollary \ref{coro4} and Proposition \ref{prop19} and take $S=\{1\}$.
\end{cor}

\begin{cor}\label{coro6}
Let P and I be graded e-weakly S-primary ideals of G-graded ring R that are not graded e-S-primary ideals of R, then $sP_e I_e=\{0\}$ for some $s\in S$.
\\
\textbf{Proof.} By Corollary \ref{coro4}, $I_e \subseteq Grad(\{0\})$, then $I_e \subseteq R_e\cap Grad(\{0\})=(Grad(\{0\}))_e$. SO, by Proposition \ref{prop19}, $sP_eI_e \subseteq sP_e(Grad(\{0\}))_e=\{0\}$ for some $s\in S$.
\end{cor}

\begin{thm}\label{thm9}
Let a graded ring homomorphism $(f)$ from $R_1$ to $R_2$, S$\subseteq$ $(R_1)_e$ be a multiplicative set of $R_1$ and $P_1, P_2$ are two graded ideal of $R_1$, $R_2$ respectively such that $P_1\cap S=\emptyset$ and $P_2\cap f(S)=\emptyset$. Then the following hold:
\\
$(i)$ If $P_1$ is a graded g-weakly S-primary ideal of $R_1$ containing  Ker(f) and  $f$  is a graded epimorphism, then $f(P_1)$ is a graded $f(g)$-weakly $f(S)$-primary ideal of $R_2$.
\\
$(ii)$ If $P_2$ is a graded $f(g)$-weakly $f(S)$-primary ideal of $R_2$ and $f$ is a graded monomorphism, then $f^{-1}(P_2)$ is a graded g-weakly S-primary ideal of $R_1$.
\\
\textbf{Proof.} $(i)$ Since $P_1$ is a graded ideal of $R_1$, so $f(P_1)$ is a graded ideal of R. Let $t\in f(P_1)\cap f(S)$ then $t=f(p)=f(s)$ when $p\in P_1 $ and $s\in S$. SO $p-s \in Ker(f) \subseteq P_1$, as a result $s\in P_1$ this is contradiction. Thus, $f(P_1)\cap f(S) = \emptyset$. For all $a,b \in (R_2)_{f(g)}$ when $g\in G$ with $0\neq ab \in f(P_1) $ . Then for some $x,y \in (R_1)_g$ such that $f(x)=a$ and $f(y)=b$, and then $0\neq f(xy)=ab \in f(P_1)$. Now, by assumption $Ker(f) \subseteq P_1$ we have got $0\neq xy \in P_1$, so $sx\in P_1$ or $sy \in Grad(P_1)$ for some $s\in S$. This implies that $f(s)a\in f(P_1)$ or $f(s)b\in Grad(P_1)$ for some  $f(s)\in f(S) \subseteq f((R_1)_e)$. Thus, $f(P_1)$ is graded f(g)-weakly f(S)-primary ideal of $R_2$.
\\
$(ii)$ Since $f^{-1}(P_2)$ is a graded ideal of $R_1$ then $f^{-1} \cap S=\emptyset$. For all $a, b \in R_g $ when $g\in G$ with $0\neq ab \in f^{-1}(P_2)$, then $0\neq f(ab)=f(a)f(b)\in P_2$ since $Ker(f)=\{0\}$, and then $f(sa)\in P_2$ or $f(sb)\in Grad(P_2))$ for some $s\in S$. Thus, $sa \in f^{-1}(P_2) $ or $sb \in Grad(f^{-1}(P_2))$. Hence $f^{-1}(P_2)$ is a graded g-weakly S-primary ideal of $R_1$.
\end{thm}

\begin{thm}\label{thm10}
  Let R be a G-graded valuation domain, then every graded ideal P of R such that $P\cap S=\emptyset$ is a graded g-weakly S-primary ideal  of R.
   \\
   \textbf{Proof.} Since P is a graded ideal of R, then  $Grad(P)= \bigcap\limits_{P\subseteq I} I$, when I denotes any graded prime ideal containing P. Since R is a G-graded valuation domain all graded prime ideals which contains P is linearly ordered, So Grad(P)=$I_\alpha$ for some graded prime ideal $I_\alpha$, and so Grad(P) is a graded prime ideal of R thus Grad(P) is a graded g-weakly S-prime ideal of R. By Proposition \ref{prop13} P is a graded weakly S-primary ideal of R.
\end{thm}.

\end{document}